\documentclass{amsart}

\usepackage{amsmath,amsthm,amssymb,amsfonts,enumerate}

\numberwithin{equation}{section} 
\theoremstyle{plain}
\newtheorem{theo+}           {Theorem}      [section]
\newtheorem{prop+}  [theo+]  {Proposition}
\newtheorem{coro+}  [theo+]  {Corollary}
\newtheorem{lemm+}  [theo+]  {Lemma}
\newtheorem{defi+}  [theo+]  {Definition}

\theoremstyle{definition}
\newtheorem{exam+}  [theo+]  {Example}
\newtheorem{rema+}  [theo+]  {Remark}

\newenvironment{theorem}{\begin{theo+}}{\end{theo+}}
\newenvironment{proposition}{\begin{prop+}}{\end{prop+}}
\newenvironment{corollary}{\begin{coro+}}{\end{coro+}}
\newenvironment{lemma}{\begin{lemm+}}{\end{lemm+}}

\begin{document}

\baselineskip 18pt
\larger[2]
\title
[Sklyanin invariant integration]
{Sklyanin invariant integration} 
\author{Hjalmar Rosengren}
\address
{Department of Mathematics\\ Chalmers University of Technology and G\"oteborg
 University\\SE-412~96 G\"oteborg, Sweden}
\email{hjalmar@math.chalmers.se}
\urladdr{http://www.math.chalmers.se/{\textasciitilde}hjalmar}

\begin{abstract}
The Sklyanin algebra admits  realizations by difference operators
acting on  theta functions. Sklyanin found an invariant
metric for the action and conjectured an explicit formula for the
corresponding reproducing kernel. We prove this conjecture, and also
 give natural biorthogonal and orthogonal bases for the representation
space. Moreover, we discuss connections with elliptic hypergeometric
series and integrals and with elliptic $6j$-symbols. 
\end{abstract}

\maketitle        

\section{Introduction}  

The study of solvable models in statistical mechanics and related
areas of physics led to the introduction of quantum groups in the
 1980's. In one of the earliest papers on the subject \cite{sk1},
Sklyanin introduced what has become known as the Sklyanin algebra.
 Its commutation relations were
obtained from the Boltzmann weights of the eight-vertex model. Since 
these involve elliptic functions, the Sklyanin algebra is an example
of an elliptic quantum group. 

The development of elliptic quantum groups has
been slow, especially regarding analytic aspects. We think here of
concrete problems in harmonic analysis, typically resulting in
explicit identities involving special functions. In a recent paper \cite{r}, we
made some progress by explaining how analytically continued 
elliptic $6j$-symbols appear in connection with the Sklyanin algebra.
These symbols  are Boltzmann weights for a 
generalization of the eight-vertex model \cite{d3}, see also \cite{d2,ft}.
In \cite{r}, they appeared as matrix elements for the change between
natural bases of finite-dimensional representations.
 The difference from the
case of Lie groups, or even simpler quantum groups, is that in that
situation  ``natural'' would mean the eigenbasis of a Lie
algebra element, whereas  for the Sklyanin algebra
one must consider a generalized eigenvalue problem $Y_1v=\lambda
Y_2v$ involving two different algebra elements.  

In the present paper we use the results of \cite{r} to further develop
 harmonic analysis on the Sklyanin algebra. 
In particular, we are interested in questions connected with invariant
integration on  a fixed representation (as opposed to questions
 connected with the Haar measure). We work with
 representations found by Sklyanin \cite{sk2}, where the algebra
 acts   by difference operators
on spaces of higher order theta functions, that is, on sections of
certain line bundles on a  torus. Sklyanin introduced a measure on
the torus which is invariant  in the sense that his algebra generators are
self-adjoint on the corresponding $\mathrm L^2$-space.

In spite of its elegance, so far there seems to have been no applications of
Sklyanin's result. At least part of the reason must be  that until now
 there has been no  example of two functions whose scalar
 product can be computed explicitly. However, for a metric to be useful it
seems necessary to have a rich class of such examples.
 In particular, one would like  to know an explicit orthogonal
basis and an explicit expression for the reproducing kernel.

Sklyanin conjectured an explicit formula for the
reproducing kernel, a problem that has  remained open. 
Our first main result, Theorem \ref{sxc}, settles this
conjecture.  In our second main result, Theorem \ref{bbt}, we give
explicit biorthogonal bases, which are ``natural'' in the sense
alluded to above.  By specialization, one may obtain orthogonal
bases, which allows us to prove another conjecture of Sklyanin,
Proposition \ref{acp},
concerning the  continuation of his scalar product in the 
 parameter $q$.
Our main tools are the alternative description of the Sklyanin algebra 
due to Rains \cite{rai} (see also the Appendix), together with
the generalized eigenvalue and tridiagonal equations from \cite{r},
expressing how Rains' operators act on our bases.

The plan of the paper is as follows. Section \ref{ps} contains
preliminaries and Section~\ref{rs} statements of the main results. The
proofs of these are given in Sections~\ref{iis}--\ref{nbs}. In the
remaining three sections  we comment briefly on
relations to other topics. In
Section~\ref{ehs} we explain how the most fundamental identity for
elliptic hypergeometric series, the Frenkel--Turaev summation, arises
naturally as the bridge between our two main results. 
In Section \ref{eis} we point out another by-product,
the evaluation of an elliptic hypergeometric double
integral, which we have not been able to reduce to known results. In
the final Section \ref{esjs}, we explain the relevance of the present work
to elliptic $6j$-symbols: it explains their self-duality or, in the
language of harmonic analysis, the duality between the spectral and
geometric variables. In the Appendix we provide some details about the
relation between Rains' and Sklyanin's difference operators.

\section{Preliminaries}
\label{ps}

\subsection{Theta functions}

Throughout,  $\tau$ and $\eta$ will be  fixed
parameters, and we let
$$p=e^{2\pi i \tau},\qquad q=e^{4\pi i\eta}.$$
We  assume that $\tau\in i\mathbb R_{>0}$, or equivalently
$0<p<1$.  The parameter $\eta$ will be either real or purely
imaginary, so that $q$ is unimodular or positive, and will later
become subject to further restrictions.

Elliptic functions are built from theta functions similarly as
rational functions are built from first degree polynomials. 
We  take as our building block Jacobi's function
$\theta_1(x|\tau)$, 
which we denote for short by $\theta(x)$. Thus,
\begin{equation}\label{tf}\begin{split}\theta(x)=\theta_1(x|\tau)&=
i\sum_{n=-\infty}^\infty(-1)^ne^{\pi i\tau(n-1/2)^2+\pi i(2n-1)x}\\
&=ip^{1/8}e^{-\pi  ix}(p,e^{2\pi ix},pe^{-2\pi
i  x};p)_\infty,\end{split}\end{equation} 
where, in general,
$$(a_1,\dots,a_n;p)_\infty=\prod_{j=0}^\infty(1-a_1p^j)\dotsm(1-a_np^j).$$
The function $\theta$  is  entire  with zeroes
$\mathbb Z+\tau\mathbb Z$.
We will often use  short-hand notation such as
\begin{subequations}\label{tsh} 
\begin{equation}\theta(a_1,\dots,a_n)=\theta(a_1)\dotsm\theta(a_n),\end{equation}
\begin{equation}\theta(a\pm b)=\theta(a+b)\theta(a-b),\end{equation}
\end{subequations}
and, for $x\in\mathbb C$, $y\in\mathbb C^n$,
$$\theta(x+\vec y)=\theta(x+y_1,\dots,x+y_n).$$

The function $\theta$ satisfies the elementary identities
$$\theta(-x)=-\theta(x), \qquad \overline{\theta(x)}=\theta(\bar x),$$
\begin{equation}\label{qp}
\theta(x+1)=-\theta(x),\qquad \theta(x+\tau)=-e^{-\pi i
(2x+\tau)}\theta(x),\end{equation}
\begin{equation}\label{dupl}
\theta(2x)=\frac{ip^{1/8}}{(p;p)_\infty^3}
\,\theta\left(x,x+\frac 12,x+\frac \tau 2,
x-\frac 12-\frac \tau 2\right),\end{equation}
 the addition formula
\begin{equation}\label{add}
\theta(x\pm y,u\pm v)=\theta(u\pm x,v\pm y)-\theta(u\pm y,v\pm x)\end{equation}
and the more general identity \cite[p.\ 451]{ww}, see also \cite[Section 4]{r2},
\begin{equation}\label{pf}
\sum_{k=1}^n\frac{\prod_{j=1}^n\theta(x_k-y_j)}{\prod_{j=1,\,j\neq
    k}^n \theta(x_k-x_j)}=0,\qquad \sum_{j=1}^n x_j=\sum_{j=1}^n y_j.
\end{equation}
By  \eqref{tf}, 
\begin{equation}\label{tdv}\theta'(0)=\lim_{x\rightarrow 0}\frac{\theta(x)}{x}=2\pi
p^{1/8}(p;p)_\infty^3.\end{equation}
We also mention the modular transformation 
\begin{equation}\label{mt}
\theta_1(x/\tau\,|-1/\tau)=-i(\tau/i)^{1/2}e^{\pi ix^2/\tau}\theta_1(x|\tau)
.\end{equation}

Another useful result is Jacobi's identity \cite[p.\ 468]{ww},  which we write
as
\begin{multline}\label{ji}
\theta\left(\vec b-\frac B2\right)=\frac12\left\{
\theta\big(\vec b\big)+\theta\left(\vec b+\frac12\right)\right.\\
\left.+e^{\pi i(\tau+B)}\theta\left(\vec
b+\frac\tau2\right)-e^{\pi i(\tau+B)}\theta\left(\vec
b+\frac12+\frac\tau2\right)
\right\},
\end{multline}
where $b=(b_1,b_2,b_3,b_4)$ and $B=b_1+b_2+b_3+b_4$. 
In algebraic geometry, \eqref{ji} is  often called
Riemann's relation \cite{m}. A direct  proof is easy; it can also be
obtained from the case $n=5$ of \eqref{pf} after substituting
$$x=\left(0,\frac12,\frac\tau2,-\frac12-\frac\tau2,-\frac B2\right),\qquad
y=\left(-b_1,-b_2,-b_3,-b_4,\frac B2\right)$$
and using \eqref{dupl} together with the identity
$$\theta\left(\frac12,\frac\tau2,-\frac12-\frac\tau2\right)=\frac{2(p;p)_\infty^3}{ip^{1/8}},$$
which follows from \eqref{dupl} after dividing with $\theta(x)$ and
letting $x\rightarrow 0$.

Occasionally, we will  write 
$$[x]=\theta(2\eta x)$$
and denote elliptic shifted factorials by
$$[x]_k=[x][x+1]\dotsm[x+k-1],$$
$$[x_1,\dots,x_n]_k=[x_1]_k\dotsm[x_n]_k.$$
Note that
$$\lim_{q\rightarrow 1}\lim_{p\rightarrow 0}\frac{[x]_k}{[y]_k}=\frac{x(x+1)\dotsm(x+k-1)}{y(y+1)\dotsm(y+k-1)},$$
which exhibits $[x]_k$ as a two-parameter deformation of the classical Pochhammer
symbol.

\subsection{Higher order theta functions}

Let $\Theta_N$ denote the space of  even theta functions of
order $2N$ with quasi-period $(1,\tau)$ 
and zero characteristics. That is,   $\Theta_N$
consists of entire functions satisfying
$$f(x+1)=f(x),\qquad f(x+\tau)=e^{-2\pi i N(2x+\tau)}f(x),\qquad f(-x)=f(x).$$
This space has dimension $N+1$ and  is
spanned by functions of the form
$$\prod_{j=1}^N\theta(a_j\pm x).$$
 In \cite{r} we were led to choose
 $\{a_j\}$ as the union of two arithmetic
progressions. More precisely, let us write
\begin{equation}\label{awb}e_k(x)=e_k^N(x;a,b)=\prod_{j=0}^{k-1}\theta(a\pm
x+2j\eta)\prod_{j=0}^{N-k-1}\theta(b\pm x+2j\eta).\end{equation}  
Then \cite[Remark 5.2]{r} $(e_k)_{k=0}^N$ form a basis for 
 $\Theta_N$ if and only if
\begin{subequations}\label{bc}
\begin{gather}a-b+2j\eta\notin\mathbb Z+\tau\mathbb Z,\qquad j=1-N,2-N,\dots,
N-1,\\
a+b+ 2j\eta\notin  \mathbb Z+\tau\mathbb Z,\qquad
j=0,1,\dots,N-1.\end{gather}
\end{subequations}
We  think of $\Theta_N$ as a deformation of the space of polynomials of
degree $\leq N$ and of $e_k(x)$ as an analogue of
$(a+x)^k(b+x)^{N-k}$. 

We note that
\begin{equation}\label{ec}
\overline{e_k^N(x;a,b)}=e_k^N(\bar x;\pm\bar a,\pm\bar
b),\end{equation}
where the plus sign is chosen for $\eta\in\mathbb R$ and the minus sign
for $\eta\in i\mathbb R$. Moreover, one has
\begin{equation}\label{eqp}\begin{split}
e_k^N(x;a+1,b)&=e_k^N(x;a,b),\\
e_k^N(x;a+\tau,b)&=e^{-2\pi ik(\tau+2a+2(k-1)\eta)}e_k^N(x;a,b),
\end{split}\end{equation}
and similarly for $b$ since $e_k^N(x;a,b)=e_{N-k}^N(x;b,a)$.

\subsection{Sklyanin algebra}

For $a=(a_1,a_2,a_3,a_4)$ such that $\sum a_i=0$, we let
 $\Delta(a)$ denote the difference operator
\begin{equation}\label{rso}
\Delta(a)f(x)=\frac{\theta(x+\vec a-\frac 12N\eta)f(x+\eta)-
\theta(x-\vec a+\frac 12N\eta)f(x-\eta)}{\theta(2x)}.\end{equation}
It is easy to check that $\Delta(a)$ preserves the space
$\Theta_N$. We note the quasi-periodicity 
\begin{equation*}\begin{split}\Delta(a_1+1,a_2-1,a_3,a_4)&
=\Delta(a_1,a_2,a_3,a_4),\\
\Delta(a_1+\tau,a_2-\tau,a_3,a_4)&=e^{2\pi i(a_2-a_1-\tau)}
\Delta(a_1,a_2,a_3,a_4),
\end{split}\end{equation*}
and similarly for the other parameters by symmetry.

The operators $\Delta(a)$ were introduced by Rains \cite{rai1}, who also
 observed \cite{rai} that they form
the degree one subspace of representations of the Sklyanin
algebra discovered by Sklyanin \cite{sk2}.  Namely, Sklyanin
introduced four operators $S_0$, $S_1$, $S_2$ and $S_3$, corresponding
to his four algebra generators. One can  check that each $S_i$ is
 proportional to an operator of the form $\Delta(a)$ and, conversely,
that any  $\Delta(a)$ is a linear combination of the $S_i$. Although we
 will not need these facts, for
the benefit of the interested reader we provide the
details in the Appendix.

In \cite[Proposition 6.2]{r}, we considered
 the action of the operators \eqref{rso} on
the  basis vectors \eqref{awb}. In the present notation, we proved that
\begin{multline}\label{gev}\Delta( a)\,e_k^N(x;a_1-\frac 12N\eta+\eta,a_2-\frac 12N\eta+\eta)\\
\begin{split}&=-\theta\left(a_1+a_2+N\eta,a_1+a_3+(2k-N)\eta,a_2+a_3+(N-2k)\eta\right)\\
&\quad\times e_k^N(x;a_1-\frac 12N\eta,a_2-\frac 12N\eta), \end{split}\end{multline}
and that
\begin{equation}\label{td}
\Delta( a)\,e_k^N(x;\lambda,\mu)=\sum_{j=k-1}^{k+1}C_j\,e_j^N(x;\lambda+\eta,\mu+\eta)
\end{equation}
for some coefficients $C_j$. We need to know   $C_{k\pm 1}$ explicitly. This can be achieved by
choosing $x=\lambda+\eta(2k-1)$ and $x=\mu+\eta(2N-2k-1)$ in
\eqref{td}, giving
\begin{equation}\label{erc}
\begin{split}
C_{k-1}&=\frac{\theta(\lambda+\vec a+\eta(2k-1-\frac 12N),2k\eta,\lambda-\mu+2k\eta)}{\theta(\lambda+\mu+2N\eta,\lambda-\mu+2(2k-N-1)\eta,\lambda-\mu+2(2k-N)\eta)},\\
C_{k+1}&=\frac{\theta(\mu+\vec a+\eta(\frac
  32N-2k-1),2(k-N)\eta,\lambda-\mu+2(k-N)\eta)}{\theta(\lambda+\mu+2N\eta,\lambda-\mu+2(2k-N)\eta,\lambda-\mu+2(2k-N+1)\eta)}.\end{split}
\end{equation}

\subsection{Involution}

Several of our results are most conveniently stated in terms of the
involution $\sigma$ on $\Theta_N$ defined by
$$(\sigma f)(x)=e^{2\pi iN(\frac 14+\frac
  \tau4+x)}f\left(x+\frac12+\frac\tau2\right)=e^{2\pi iN(\frac 14+\frac
  \tau4-x)}f\left(x-\frac12-\frac\tau2\right). $$
It is easy to check that $\sigma$ preserves $\Theta_N$ and that
$\sigma\circ\sigma=\operatorname{id}$.  

We  mention the easily verified identity
\begin{equation}\label{se}\begin{split}\sigma e_k^N(x;a,b)&=e^{2\pi i(ak-b(N-k)+(N-1)(2k-N)\eta+\frac 14
  N(\tau-1))}\\
&\quad\times e_k^N(x;a+\frac 12+\frac\tau 2,b-\frac
12-\frac\tau 2). \end{split}\end{equation}

We also need to know how  $\Delta(a)$ behaves under
conjugation by $\sigma$. One first computes
\begin{multline}\label{psd}(\sigma\Delta(a)\sigma f)(x)=\frac{e^{\pi i(\tau+4x)}}{\theta(2x)}\left\{e^{-2\pi iN\eta}\theta(x+\vec a+\frac12+\frac\tau2-\frac
  12N\eta)f(x+\eta)\right.\\
\left.-e^{2\pi iN\eta}
\theta(x-\vec a+\frac12+\frac\tau2+\frac 12N\eta)f(x-\eta)\right\},
 \end{multline}
which may be rewritten in the form
\begin{multline}\label{sd}
\sigma\circ\Delta(a)\circ\sigma\\
=e^{2\pi i(a_1+a_3+\frac
  \tau2)}\Delta\left(a_1+\frac 12+\frac\tau 2,a_2+\frac 12-\frac\tau
2, a_3-\frac 12+\frac\tau 2, a_4-\frac
12-\frac\tau 2\right).
\end{multline}
The apparent loss of symmetry  is needed
 to preserve the condition $\sum a_i=0$.

\subsection{Invariant integration}

The following metric on $\Theta_N$ was introduced by Sklyanin \cite{sk2}:
\begin{equation}\label{im}
\langle f,g\rangle=\iint_{\mathbb C/(\mathbb Z+\tau\mathbb Z)}
f(u)\,\overline{g(u)}\,M(u,\bar u)
\,dxdy,\end{equation}
where $u=x+iy$ and
$$M(u,v)=\frac{\theta(2u)\theta(2v)}{e^{2\pi
    iu(N+2)}\prod_{k=0}^{N+1}\theta(u\pm v+(2k-N-1)\eta+\frac
    12+\frac\tau 2)}.$$
It can be viewed as a deformation of the standard 
$\mathrm{SU}(2)$-invariant metric on polynomials of degree $\leq N$:
$$\text{Const}
\iint_{\mathbb C}\frac{f(u)\,\overline{g(u)}}{(1+|u|^2)^{N+2}}\,dxdy. $$

It is easy to check that,  for 
 $f,g\in\Theta_N$, the integrand in \eqref{im}
has indeed double-period $(1,\tau)$. 
The weight $M(u,\bar u)$ is free from poles if
\begin{equation}\label{wa}
(2k-N-1)\eta\notin\left(\mathbb Z+\frac 12+i\mathbb R\right)\cup
\left(\tau(\mathbb Z+\frac 12)+\mathbb R\right),\qquad
k=0,1,\dots,N+1,\end{equation}
and  it is non-negative if $\eta\in\mathbb R\cup i\mathbb R$, since
then
$$M(u,\bar u)=\frac{p^{\frac14(N+2)}}{(p;p)_\infty^{2N+2}}\,
\left|\theta(2u)\right|^2\prod_{k=0}^{N+1}\left|\frac{1}{(-\sqrt
    p\,e^{2\pi i(u\pm\bar u+(2k-N-1)\eta)};p)_\infty}\right|^2.
$$
 When
both these conditions are satisfied, we have a genuine scalar
product. In particular, this happens for 
$\eta\in \mathbb R_N\cup \mathbb I_N$, where 
$$\mathbb R_N=\{\eta\in\mathbb R;|\eta|<1/2(N+1)\},\qquad 
\mathbb I_N=\{\eta\in i\mathbb R;|\eta|<\tau/2(N+1)i\}.$$

\section{Statement of results}
\label{rs}

Our main tool is the
invariance of the metric \eqref{im}
with respect to the operators \eqref{rso}. 
Sklyanin proved that for $\eta\in\mathbb R_N$ 
his operators $S_i$  are self-adjoint. 
Working with the more general operators \eqref{rso} allows one
to simplify the proof, since one need not then
consider the  $S_i$ one by one. Moreover, we want to extend Sklyanin's
result to the case  $\eta\in
\mathbb I_N$. Although the two cases are related 
by the modular transformation \eqref{mt}, we prefer to treat
them in parallel.  For these two reasons,
and since Sklyanin's presentation is rather sketchy, we 
give a detailed proof in
 Section \ref{iis}.

\begin{proposition}\label{iil}
Assume that  $\eta\in\mathbb R_N$, and let
 $\Delta(a)$, $\sum a_i=0$, be an operator of the form
\eqref{rso}. Then its adjoint with respect to the metric \eqref{im}
is given by 
$\Delta(a)^\ast=-\sigma\circ\Delta(-\bar a)\circ\sigma$. Similarly, if 
 $\eta\in
\mathbb I_N$,  then  $\Delta(a)^\ast=\sigma\circ\Delta(\bar a)\circ\sigma$.
\end{proposition}

It is easy to check that $\sigma^\ast=\sigma$, which shows that
Proposition \ref{iil}
 is consistent with $\Delta(a)^{\ast\ast}=\Delta(a)$. Note
also that, by \eqref{sd}, there exist in both cases a constant $C$ and
parameters $b_i$ such 
that $\Delta(a)^\ast=C\Delta(b)$. However, the formulation involving $\sigma$
is   more convenient for our purposes. 

Using  \eqref{sg} and \eqref{sd}, one may check that Proposition
\ref{iil} agrees with the result of Sklyanin mentioned above.

\begin{corollary}
For  $\eta\in\mathbb R_N$,  one has
 $S_i^\ast=S_i$, $i=0,1,2,3$, whereas for  $\eta\in
\mathbb I_N$,   $S_i^\ast=-S_i$, $i=0,1,2,3$.
\end{corollary}

Our first main result concerns  the reproducing kernel
 for $\Theta_N$ with respect to the metric \eqref{im}.
 Sklyanin conjectured that it is given, up to a
multiplicative constant, by
$$K_v(u)=K(u,\bar v)=e^{2\pi i uN}\prod_{k=0}^{N-1}
\theta\left(u\pm\bar v+(2k-N+1)\eta+\frac 12+\frac\tau 2\right).$$ 
Note that
\begin{equation}\label{ksymm}K(u,\bar v)=\overline{K(v,\bar
    u)}=K(\bar v,u).\end{equation}
In Section \ref{scs} we will prove Sklyanin's conjecture, and in Section
\ref{ccs} the multiplicative constant  will be computed. 
 We summarize the result as follows.

\begin{theorem}\label{sxc}
For $\eta\in\mathbb R_N\cup\mathbb I_N$, the reproducing kernel of 
$\Theta_N$ with the metric \eqref{im} is given by $C^{-1}K_v$, where
\begin{equation}\label{c} C=\frac{2\eta p^{3/8}}{\theta(2(N+1)\eta)(p;p)_\infty^3}.\end{equation}
\end{theorem}

When $\eta=0$, the expression for $C$ has a removable singularity and
should be interpreted as, cf.\   \eqref{tdv},
$$\lim_{\eta\rightarrow 0}C=\frac{p^{1/4}}{2\pi(N+1)(p;p)_\infty^6}.$$

We  have also obtained biorthogonal bases for $\Theta_N$. Indeed, we can
find the  dual of any basis of the form \eqref{awb}.

\begin{theorem}\label{bbt}
Suppose that 
$$e_k(x)=e_k^N(x;a_1-\frac 12N\eta,a_2-\frac 12N\eta),\qquad 
k=0,\dots,N,$$
 form a
basis for the space $\Theta_N$. Let 
$$
f_k(x)=\sigma e_k^N(x;\mp \bar a_2-\frac 12N\eta+\eta,\mp \bar a_1-\frac 12N\eta+\eta),$$
 where the minus sign is chosen if $\eta\in\mathbb R_N$ and the
plus sign if $\eta\in\mathbb I_N$. Then
$$\langle e_k,f_l\rangle=C\Gamma_k\,\delta_{kl},$$
where   $C$ is given by \eqref{c} and
$$\Gamma_k=e^{\pi iN(\tau-1)/2}\frac{[\lambda]}{[\lambda+2k]}\frac{[1,\lambda+N+1]_k}{[-N,\lambda]_k}\,[\lambda+1,(a_1+a_2-N\eta)/2\eta]_N,$$
with
$\lambda=(a_1-a_2-2N\eta)/2\eta$. 
\end{theorem}

In the case $\eta=0$, the expression for $\Gamma_k$ should be
interpreted as the limit
$$\lim_{\eta\rightarrow 0}\Gamma_k=e^{\pi iN(\tau-1)/2}
\frac{(-1)^k}{\binom{N}{k}}\,\theta(a_1+ a_2)^N\theta(a_1-a_2)^N.
 $$
We also observe that, by \eqref{ec}, for $\eta\in\mathbb R_N$ and
$\eta\in\mathbb I_N$ alike we have
\begin{equation}\label{fc}
\overline{f_k(x)}=e^{2\pi iN(\bar x-\frac14+\frac\tau4)}e_k^N(\bar x+\frac12+\frac\tau2;-a_2-\frac12N\eta+\eta,-a_1-\frac12N\eta+\eta).
\end{equation}

The 
 proof of  Theorem \ref{bbt} is divided into two parts.
For the biorthogonality we use the generalized eigenvalue equation
 \eqref{gev} and for the norm computation the generalized tridiagonal
 equation \eqref{td}. The details are given in
Sections \ref{bbs} and \ref{nbs}, respectively. 

Sklyanin conjectured \cite[p.\ 277]{sk2} that his 
 metric extends from
 $\eta\in\mathbb R_{N}$ to the larger parameter range 
$\eta\in\mathbb R_{N-1}$. This follows quite easily from Theorem
 \ref{bbt}.

\begin{proposition}\label{acp}
The analytic continuation in $\eta$ of
 $C^{-1}\langle f,g\rangle$ is positive definite for
$\eta\in\mathbb 
R_{N-1}\cup\mathbb I_{N-1}$.  
\end{proposition}

\begin{proof}
We first specialize the parameters in Theorem \ref{bbt} to
obtain an orthogonal basis.
By \eqref{se},  $f_k(x)$ is proportional to 
$$e_k^N(x;\mp \bar a_2-\frac 12N\eta+\eta+\frac12+\frac\tau2,\mp
\bar a_1-\frac 12N\eta+\eta-\frac12-\frac\tau2).$$ 
In particular, if   $a_2=\mp \bar
a_1+\eta-\frac12-\frac\tau2$ then, using also \eqref{eqp},
$e_k$ and $f_k$ are proportional, so that $(e_k)_{k=0}^N$ is 
orthogonal. We know that  $C^{-1}\|e_k\|^2$ is
given by an exponential factor times $\Gamma_k$, and that it is
 positive for  $\eta\in\mathbb R_{N}\cup\mathbb
I_{N}$. By continuity, for generic $a_1$ it must
remain positive as long as $[1]_k/[-N]_k$ is well-defined and
non-zero. This is indeed true for  $\eta\in\mathbb R_{N-1}\cup\mathbb
I_{N-1}$. 
\end{proof}

\section{Proof of Proposition \ref{iil}}
\label{iis}

Consider first the case  $\eta\in\mathbb R_N$. We write
\begin{equation*}\begin{split}\langle
    \Delta(a)f,g\rangle&=\iint_{v=\bar u}
\frac{\theta(u+\vec a-\frac 12N\eta)f(u+\eta)-
\theta(u-\vec a+\frac 12N\eta)f(u-\eta)}{\theta(2u)}\\
&\quad\times\overline{g(\bar v)}\,M(u,v)\,dxdy\\
&=\iint_{v=\bar u-\eta}
\frac{\theta(u+\vec a-\frac 12(N+2)\eta)f(u)
\overline{g(\bar v)}\,M(u-\eta,v)}{\theta(2u-2\eta)}\,\,dxdy\\
&\quad - \iint_{v=\bar u+\eta}
\frac{\theta(u-\vec a+\frac 12(N+2)\eta)f(u)
\overline{g(\bar v)}\,M(u+\eta,v)}{\theta(2u+2\eta)}\,\,dxdy.
\end{split}\end{equation*}
 We wish to replace the contours of integration with
 $v= \bar u$. In general, we have that
$$\iint_{v=\bar u-\gamma}
\frac{\theta(u+\vec a-\frac 12(N+2)\eta)f(u)
\overline{g(\bar v)}\,M(u-\eta,v)}{\theta(2u-2\eta)}\,\,dxdy$$
is independent of $\gamma$ as long as we avoid the poles of
$M(u-\eta,\bar u-\gamma)$, that is, for
$$\gamma\notin \mathbb R+\tau\left(\mathbb Z+\frac 12\right),$$
$$\gamma-\eta+(2k-N-1)\eta\notin \mathbb Z+\frac 12+i\mathbb R, \qquad
k=0,1,\dots,N+1.$$
In particular, the region of analyticity containing
$\gamma=\eta$ is given by
$$|\operatorname{Im}(\gamma)|<\frac{\tau}{2i},\qquad 
|\operatorname{Re}(\gamma)-\eta\pm(N+1)\eta|<1/2.$$ 
If we make the temporary assumption $\eta\in\mathbb R_{N+1}$,
 then this region contains $\gamma=0$, so
that the contour  may indeed be replaced with $v=\bar
u$. This also holds for the integral on $v=\bar u+\eta$, giving
\begin{multline}\label{das}\langle
    \Delta(a)f,g\rangle=\iint_{v=\bar u}f(u)\overline{g(\bar v)}
\left\{\frac{\theta(u+\vec a-\frac
    12(N+2)\eta)M(u-\eta,v)}{\theta(2u-2\eta)}\right.\\
\quad\left.
-\frac{\theta(u-\vec a+\frac 12(N+2)\eta)M(u+\eta,v)}{\theta(2u+2\eta)}
\right\}\,dxdy.\end{multline}

Using \eqref{psd}, the same argument, still assuming $\eta\in\mathbb R_{N+1}$, gives
\begin{multline}\label{dbs}\langle f,\sigma\Delta(\bar b)\sigma g\rangle\\
=
\iint_{v=\bar u}f(u)\overline{g(\bar v)}\,e^{\pi i(\tau-4v)}
\left\{e^{2\pi i(N+2)\eta}\frac{\theta(v+\vec b+\frac12-\frac\tau2-\frac
    12(N+2)\eta)M(u,v-\eta)}{\theta(2v-2\eta)}\right.\\
\left.-e^{-2\pi i(N+2)\eta}
\frac{\theta(v-\vec b+\frac12-\frac\tau2+\frac 12(N+2)\eta)M(u,v+\eta)}{\theta(2v+2\eta)}
\right\}\,dxdy.\end{multline}
We are thus reduced to proving
\begin{multline}\label{pmrec}\frac{\theta(u+\vec a-\frac
    12(N+2)\eta)M(u-\eta,v)}{\theta(2u-2\eta)}-\frac{\theta(u-\vec a+\frac 12(N+2)\eta)M(u+\eta,v)}{\theta(2u+2\eta)}\\
=e^{i\pi(\tau-4v)}
\left\{e^{-2\pi i(N+2)\eta}\frac{\theta(v+\vec a+\frac12-\frac\tau2+\frac
    12(N+2)\eta)M(u,v+\eta)}{\theta(2v+2\eta)}\right.\\
\left.-e^{2\pi
    i(N+2)\eta}\frac{\theta(v-\vec a+\frac12-\frac\tau2-\frac 12(N+2)\eta)M(u,v-\eta)}{\theta(2v-2\eta)}\right\}.
\end{multline}

Multiplying \eqref{pmrec} with $e^{2\pi
  i(N+2)u}\prod_{k=0}^{N+2}\theta(u\pm v+(2k-N-2)\eta+\frac
  12+\frac\tau 2)$ and simplifying gives
 \begin{multline}\label{mrec}
\theta(2v)\Big\{e^{2\pi i(N+2)\eta}
\,\theta(\textstyle u\pm v+(N+2)\eta+\frac 12+\frac\tau 2,u+\vec a-\frac 12(N+2)\eta)
\\
-e^{-2\pi i(N+2)\eta}
\,\theta(\textstyle u\pm v-(N+2)\eta+\frac 12+\frac\tau 2,u-\vec
a+\frac 12(N+2)\eta)\Big\}
\\
\begin{split}&=e^{i\pi(\tau-4v)}\,
\theta(2u)\\
&\quad\times\Big\{e^{-2\pi i(N+2)\eta}
\theta(\textstyle
u\pm (v-(N+2)\eta)+\frac 12+\frac\tau 2,v+\vec a+\frac 12-\frac\tau 2+\frac 12(N+2)\eta)\end{split}\\
-e^{2\pi i(N+2)\eta}
\theta(\textstyle u\pm (v+(N+2)\eta)+\frac 12+\frac\tau
  2,v-\vec a+\frac 12-\frac\tau 2-\frac 12 (N+2)\eta) \Big\}.
\end{multline}
This is equivalent to the case  $n=4$ of \eqref{pf}, which we rewrite  as
\begin{multline*}
\theta(x_4-x_3)\left\{\theta(x_2-x_3,x_2-x_4,x_1-\vec y)-
\theta(x_3-x_1,x_4-x_1,x_2-\vec y) \right\}\\
=\theta(x_2-x_1)\left\{\theta(x_4-x_1,x_2-x_4,x_3-\vec y)
-\theta(x_3-x_1,x_2-x_3,x_4-\vec y) \right\},
\end{multline*}
valid for $y=(y_1,y_2,y_3,y_4)$ with  $\sum_ix_i=\sum_iy_i$.
Indeed, substituting
\begin{multline*}(x_1,x_2,x_3,x_4)
 =\left(-u+\frac 12{(N+2)\eta},u+\frac 12{(N+2)\eta},\right.\\
 \left. -v-\frac 12+\frac\tau
 2-\frac12 {(N+2)\eta},v+\frac 12-\frac\tau 2-\frac12 {(N+2)\eta}\right), 
\end{multline*}
letting $y_i=a_i$ and repeatedly using \eqref{qp}, one obtains  \eqref{mrec}.

We have assumed  $\eta\in\mathbb R_{N+1}$,
 but the result extends to $\eta\in\mathbb R_{N}$. Namely, 
$\langle\Delta(a)f,g\rangle+\langle f,\sigma\Delta(-\bar
a)\sigma g\rangle$ is analytic in $\eta$ as long as \eqref{wa} is
satisfied; thus, it is zero for $\eta\in\mathbb R_N$.

Repeating the calculation for imaginary $\eta$, one sees that 
\eqref{das} holds also for  $\eta\in \mathbb I_{N+1}$, while in
\eqref{dbs}, $\eta$ should be replaced by $\bar\eta=-\eta$ everywhere
on the right-hand side. However, this has the same effect as replacing
$b$ by $-b$ and multiplying the whole expression with $-1$. Thus, we are
again reduced to the identity \eqref{pmrec}. The extension from
$\eta\in \mathbb I_{N+1}$ to $\eta\in \mathbb I_{N}$ follows as before.

\section{Proof of Sklyanin's conjecture}
\label{scs}

In this section we prove Sklyanin's conjecture, that is, that
\begin{equation}\label{rki}f(u)=C^{-1}\langle f,K_u\rangle,\qquad
  f\in\Theta_N,\quad u\in\mathbb C\end{equation}
for some constant $C$.

We  need to know that the kernels $(K_u)_{u\in\mathbb C}$
span $\Theta_N$. This can be seen, for instance, by considering the
elements
$$e_k^N(x;a,2\eta-a)=\prod_{j=k-N}^{k-1}\theta(a\pm x+2j\eta),\qquad
k=0,1,\dots, N.$$ 
These are all proportional to a kernel $K_u(x)$. Moreover, by \eqref{bc}, for
generic $a$ they form a basis for $\Theta_N$ as long as
$2j\eta\notin\mathbb Z+ \tau\mathbb Z$, $1\leq j\leq N$; in
particular, this happens  for $\eta\in\mathbb R_N\cup\mathbb I_N$ (and even for 
$\eta\in\mathbb R_{N-1}\cup\mathbb I_{N-1}$), $\eta\neq 0$. 

It is thus enough to prove \eqref{rki} for $f=K_v$, that is, introducing the kernel
$$\Phi_v(u)=\Phi(u,\bar v)=\langle K_v,K_u\rangle,$$
that
\begin{equation}\label{kj}\Phi_v=CK_v.\end{equation} 

Our main tool for proving \eqref{kj} is the existence of Sklyanin
algebra elements that act nicely on the kernel $K_v$. Note first that,
by \eqref{ksymm}, we may write
$$K_v(x)=e^{2\pi i\bar vN}e_N^N(x;a_1+\eta-\frac 12N\eta,
a_2+\eta-\frac 12N\eta),$$
where  $a_1=\bar v-\frac 12 N\eta+\frac12+\frac\tau2$ and $a_2$
is arbitrary. Thus, \eqref{gev} gives
$$\Delta(a)K_v=-e^{2\pi
  iN\eta}\theta(a_1+a_2+N\eta,a_1+a_3+N\eta,a_2+a_3-N\eta)\, K_{v-\bar\eta},$$
where  $a_3$ is arbitrary and $a_4=-a_1-a_2-a_3$.

Assume first that $\eta\in\mathbb R_N$.
We consider the equality
\begin{equation}\label{ss}
\langle \Delta(a)K_v,K_u\rangle=\langle K_v,\Delta(a)^\ast K_u\rangle,
\end{equation}
where we know that the left-hand side is
\begin{equation}\label{l}-e^{2\pi
  iN\eta}\theta(a_1+a_2+N\eta,a_1+a_3+N\eta,a_2+a_3-N\eta)
\,\Phi(u,\bar v-\eta).\end{equation}

As for the right-hand side, we use  Proposition
  \ref{iil} together with the identity
$$\sigma K_u(x)=e^{\pi iN(1-\tau)/2}e_0^N(x;w,\bar u+(1-N)\eta),$$
with $w$ arbitrary. Choosing $a_2=\frac 12N\eta-u$, we find that the
right-hand side of \eqref{ss} equals
\begin{equation}\label{r}-\theta(a_1+a_2-N\eta,a_1+a_3+N\eta,a_2+a_3-N\eta)\,\Phi(u-\eta,\bar v).\end{equation}
Identifying \eqref{l} and \eqref{r} we obtain  the
difference equation
\begin{equation}\label{phr}\Phi(u,\bar v-\eta)=
e^{2\pi iN\eta}
\frac{\theta(u-\bar v+N\eta+\frac12+\frac\tau2)}{\theta(u-\bar v-N\eta+\frac12+\frac\tau2)}\,\Phi(u-\eta,\bar v),\end{equation}
which, after replacing $\bar v$ with $\bar v+\eta$ and iterating yields
\begin{equation}\label{ir}
\Phi(u,\bar v)=e^{2\pi iN\eta k}\prod_{j=1}^{k}\frac{\theta(u-\bar
  v-(2j-1-N)\eta+\frac12+\frac\tau2)}{\theta(u-\bar
  v-(2j-1+N)\eta+\frac12+\frac\tau2)}\, \Phi(u-k\eta,\bar v+k\eta).
\end{equation}

We now plug $u=u_k=\bar v+(2k-1-N)\eta-\frac12-\frac\tau2$ into \eqref{ir},
which makes the numerator zero. The denominator is then
$\prod_{j=1}^k\theta(2(k-j-N)\eta)$, which is non-zero if $1\leq k\leq
N$ and $\eta\in\mathbb R_N$, $\eta\neq 0$. Thus, 
 $\{u_k\}_{k=1}^N$ are zeroes of the 
function $\Phi_v$. Since
$\Phi_v\in\Theta_N$, the points $\pm u_k+\mathbb Z+\tau\mathbb Z$ are
also zeroes. These are precisely the zeroes of $K_v$, so $\Phi_v/K_v$
is a $(1,\tau)$-periodic entire function, and thus constant by
Liouville's theorem.

In the case $\eta\in\mathbb I_N$, $\eta\neq 0$, the same proof goes through,
although one then finds instead of \eqref{phr} the equation
$$\Phi(u,\bar v-\eta)=
e^{-2\pi iN\eta}
\frac{\theta(u+\bar v-N\eta+\frac12+\frac\tau2)}{\theta(u+\bar v+N\eta+\frac12+\frac\tau2)}\,\Phi(u+\eta,\bar v).$$

For the above proof to work it is essential that  $\eta\neq 0$,
though the case $\eta=0$ is included by continuity, cf.\ the remark
following Theorem \ref{sxc}.

\section{Computation of the constant}
\label{ccs}

Knowing that the constant $C$ of Theorem \ref{sxc} exists, we shall now
compute it. To this end, let $(e_k)_{k=0}^N$ be an orthonormal basis of
the space $\Theta_N$, so that
$$\frac 1C\,K(u,\bar v)=\sum_{k=0}^Ne_k(u)\,\overline{e_k(v)}. $$
If we put $u=v$ in this identity and integrate we obtain
$$\frac 1C\iint K(u,\bar u)M(u,\bar u)\,dxdy=\sum_{k=0}^N\|e_k\|^2=N+1,$$
that is,
\begin{equation*}\begin{split}C&=\frac {1}{N+1}\iint K(u,\bar u)M(u,\bar
u)\,dxdy\\
&=\frac {e^{-4\pi i(N+1)\eta}}{N+1}\iint\frac{\theta(2u)\theta(2\bar
  u)}{\theta(u\pm
  \bar u\pm\gamma)}\,dxdy,\end{split}\end{equation*}
where $\gamma=(N+1)\eta+\frac 12+\frac\tau 2$. Note that
$\eta\in\mathbb R_N\cup\mathbb I_N$ means that
$0<\operatorname{Re}(\gamma)<1$, $0<\operatorname{Im}(\gamma)<\tau/i$.

Applying  the following lemma now completes the proof of
Theorem \ref{sxc}.

\begin{lemma}\label{fil}
Suppose $0<\operatorname{Re}(\gamma)<1$ and
$0<\operatorname{Im}(\gamma)<\tau/i$.  Then 
$$\iint\frac{\theta(2u)\theta(2\bar u)}{\theta(u\pm
  \bar u\pm\gamma)}\,dxdy=\frac{
  p^{-1/8}(2\gamma-1-\tau)}{\theta(2\gamma)(p;p)_\infty^3}.$$ 
\end{lemma}

Our proof will be similar to the proof of
 \cite[Lemma 3.3]{rai1}.
We first      note that 
 applying $\frac{\partial}{\partial x}\big|_{x=y}$ to both sides of
\eqref{add} gives
\begin{multline*}\frac{\theta(u\pm v)}{\theta(u\pm y)\theta(v\pm
  y)}\\
=\frac{1}{\theta'(0)\theta(2y)}\left(\frac{\theta'(u+y)}{\theta(u+y)}-\frac{\theta'(u-y)}{\theta(u-y)}+\frac{\theta'(v-y)}{\theta(v-y)}-\frac{\theta'(v+y)}{\theta(v+y)} \right).\end{multline*}
Substituting   $(u,v,y)\mapsto(2x,2iy,\gamma)$ 
and using 
\eqref{tdv}, we find that  our integrand equals
$$\frac{p^{-1/8}}{2\pi\theta(2\gamma)(p;p)_\infty^3}\left(\frac{\theta'(2x+\gamma)}{\theta(2x+\gamma)}-\frac{\theta'(2x-\gamma)}{\theta(2x-\gamma)}+\frac{\theta'(2iy-\gamma)}{\theta(2iy-\gamma)}-\frac{\theta'(2iy+\gamma)}{\theta(2iy+\gamma)}
\right). $$
We are thus reduced to computing integrals of the form
$$I(\gamma)=\int_0^1\frac{\theta'(2x+\gamma)}{\theta(2x+\gamma)}\,dx,\qquad
\gamma\notin\mathbb R+\tau\mathbb Z,$$
$$
J(\gamma)=\int_{0}^{\tau/i}\frac{\theta'(2iy+\gamma)}{\theta(2iy+\gamma)}\,dy,
\qquad\gamma\notin\mathbb  Z+i\mathbb R.  
$$

\begin{lemma}\label{ijl}
For $0<\operatorname{Im}(\gamma)<\tau/i$ one has
$I(\gamma)=-i\pi$. For other values of $\gamma$, $I(\gamma)$ is determined by
 $I(\gamma+\tau)=-2\pi
i+I(\gamma)$. Similarly, for $0<\operatorname{Re}(\gamma)<1$
one has $J(\gamma)=2\pi(\frac 12-\gamma-\tau)$; 
for other values, $J(\gamma)$ is determined by
$J(\gamma)=J(\gamma+1)$.
\end{lemma}

Although Lemma \ref{ijl} is easily deduced from known results, we include a
proof for completeness. 
The functional equations for $I$ and $J$ follow from \eqref{qp}. By
analyticity, it is then enough to assume 
$\operatorname{Im}(\gamma)=\tau/2i$ and
$\operatorname{Re}(\gamma)=1/2$, respectively. 

In the case 
$\operatorname{Im}(\gamma)=\tau/2i$, we may  write
$$\theta(2x+\gamma)=i p^{1/8}(p;p)_\infty e^{-\pi
  i(2x+\gamma)}\left|(e^{2\pi i(2x+\gamma)};p)_\infty\right|^2, $$
and thus
$$\frac{\theta'(2x+\gamma)}{\theta(2x+\gamma)}=-i \pi +\frac{d}{dx}\,\log\left|(e^{2\pi i(2x+\gamma)};p)_\infty\right|.
$$
Since the last term  is $1$-periodic, we obtain indeed
$I(\gamma)=-i\pi$. 

Similarly, if $\operatorname{Re}(\gamma)=1/2$ we may write
$$\theta(2iy+\gamma)=ip^{1/8}e^{-\pi i(2iy+\gamma)}(p,-e^{-2\pi(2y+\delta)},-pe^{2\pi(2y+\delta)};p)_\infty,$$
where $\delta=(\gamma-1/2)/i\in\mathbb R$. This gives
$$\frac{\theta'(2iy+\gamma)}{\theta(2iy+\gamma)}=-\pi i+\frac 1{2i}\frac{d}{dy}\,\log(-e^{-2\pi(2y+\delta)},-pe^{2\pi(2y+\delta)};p)_\infty,$$
which is integrated to
\begin{equation*}\begin{split}
J(\gamma)&=\int_{0}^{\tau/i}\frac{\theta'(2iy+\gamma)}{\theta(2iy+\gamma)}\,dy=-\pi\tau+\frac
1{2i}\,\log\frac{(-p^2e^{-2\pi\delta},-p^{-1}e^{2\pi\delta};p)_\infty}{(-e^{-2\pi\delta},-pe^{2\pi\delta};p)_\infty}\\
&=2\pi\left(\frac 12-\gamma-\tau\right).\end{split}\end{equation*}

Using Lemma \ref{ijl}, we can
 now compute the integral in Lemma \ref{fil} as
$$\frac{p^{-1/8}}{2\pi\theta(2\gamma)(p;p)_\infty^3}\left(\frac\tau i(I(\gamma)-I(-\gamma))+J(-\gamma)-J(\gamma)\right),
$$
where $I(\pm\gamma)=\mp \pi i$, $J(\pm\gamma)=2\pi(\pm(\frac
12-\gamma)-\tau)$, 
which gives the desired result after simplification.

\section{Biorthogonality}\label{bbs}

In this section we begin the proof of  Theorem \ref{bbt} by showing that
 $\langle e_k,f_l\rangle=0$ for
$k\neq l$. Recall that
$$e_k(x)=e_k^N(x;a_1-\frac 12N\eta,a_2-\frac 12N\eta),$$
$$f_k(x)=\sigma e_k^N(x;\mp \bar
a_2-\frac 12N\eta+\eta,\mp \bar
a_1-\frac 12N\eta+\eta),$$
where the minus sign is chosen for $\eta\in\mathbb R_N$ and the plus sign for
$\eta\in\mathbb I_N$. We will also write
$$e_k^+(x)=e_k^N(x;a_1-\frac 12N\eta+\eta,a_2-\frac 12N\eta+\eta),$$
$$f_k^-(x)=\sigma e_k^N(x;\mp \bar
a_2-\frac 12N\eta,\mp \bar
a_1-\frac 12N\eta).$$

Consider the identity
$$\langle \Delta(a)e_k^+,f_l\rangle=\mp\langle e_k^+,\sigma\Delta(\mp
\bar a)\sigma f_l\rangle,$$
where $a_3$  is  arbitrary  and  $a_4$ is fixed by $\sum
a_i=0$. By \eqref{gev}, the left-hand side is given by
$$-\theta(a_1+a_2+N\eta,a_1+a_3+(2k-N)\eta,a_2+a_3+(N-2k)\eta)\,\langle
e_k,f_l\rangle$$ 
and the right-hand side by
$$-\theta(a_1+a_2-N\eta,a_1+a_3+(2l-N)\eta,a_2+a_3+(N-2l)\eta)\,\langle
e_k^+,f_l^-\rangle.$$
Assuming that the denominator is non-zero, this gives
\begin{multline}\label{ns}
\langle
e_k,f_l\rangle\\
=\frac{\theta(a_1+a_2-N\eta,a_1+a_3+(2l-N)\eta,a_2+a_3+(N-2l)\eta)}{\theta(a_1+a_2+N\eta,a_1+a_3+(2k-N)\eta,a_2+a_3+(N-2k)\eta)}
\,\langle
e_k^+,f_l^-\rangle.
\end{multline}

We  now choose 
$a_3=(2l-N)\eta-a_2$. This gives $\langle e_k,f_l\rangle=0$, as long
as  the denominator in \eqref{ns} is non-zero, that is, for
$$a_1+a_2+N\eta,\ a_1-a_2+2(k+l-N)\eta,\ 2(l-k)\eta\notin\mathbb
Z+\tau\mathbb Z. $$
 If  $k\neq l$,  the last condition holds for $\eta\in\mathbb
 R_N\cup\mathbb I_N$, $\eta\neq 0$.  
The other two conditions follow from the assumption that
$(e_k)_{k=0}^N$ form a basis, that is, that
 \eqref{bc} holds with $(a,b)$ replaced by $(a_1-\frac
12N\eta,a_2-\frac 12N\eta)$.

This proves that $\langle e_k,f_l\rangle=0$ when $k\neq l$. As a
by-product, choosing $k=l$ in \eqref{ns} gives the identity
\begin{equation}\label{nns}\langle
e_k,f_k\rangle\\
=\frac{\theta(a_1+a_2-N\eta)}{\theta(a_1+a_2+N\eta)}
\,\langle
e_k^+,f_k^-\rangle,\end{equation}
which will be used in the next section.

\section{Norm computation}\label{nbs}

In this section we complete the proof of Theorem \ref{bbt} by 
computing the constants $\Gamma_k=C^{-1}\langle e_k,f_k\rangle$. To
this end, we  consider the equality  
\begin{equation}\label{asp}\langle \Delta(A)e_k,f_{k+1}^-\rangle=\langle e_k,\Delta(A)^\ast f_{k+1}^-\rangle, \end{equation}
where  $A=(A_1,A_2,A_3,A_4)$ is arbitrary.
 By \eqref{td} and Proposition \ref{iil}, 
$$\Delta(A)e_k=c_{k-1}e_{k-1}^+ + c_{k}e_{k}^+ + c_{k+1}e_{k+1}^+,$$
 $$\Delta(A)^\ast f_{k+1}^-=d_{k}f_{k} + d_{k+1}f_{k+1} + d_{k+2}f_{k+2},$$
for some constants $c_i$ and $d_i$. By the results of the previous
section, \eqref{asp} then reduces to
$$c_{k+1}\langle e_{k+1}^+,f_{k+1}^-\rangle=\bar d_{k}\langle e_k,f_{k}\rangle,  $$
which, by \eqref{nns}, yields the recursion
$$\Gamma_{k+1}=\frac{\theta(a_1+a_2-N\eta)}{\theta(a_1+a_2+N\eta)}\frac{\bar
  d_{k}}{c_{k+1}}\,\Gamma_k. $$

Using \eqref{erc}, we compute
$$c_{k+1}=\frac{\theta(a_2+\vec
  A+\eta(N-2k-1),2(k-N)\eta,a_1-a_2+2(k-N)\eta)}
{\theta(a_1+a_2+N\eta,a_1-a_2+2(2k-N)\eta,a_1-a_2+2(2k-N+1)\eta)}, $$
$$\bar d_k=\frac{\theta(a_2+\vec
  A+\eta(N-2k-1),2(k+1)\eta,a_1-a_2+2(k+1)\eta)}
{\theta(a_1+a_2-N\eta,a_1-a_2+2(2k-N+1)\eta,a_1-a_2+2(2k-N+2)\eta)},$$ 
 giving
$$\Gamma_{k+1} =
\frac{\theta(2(k+1)\eta,a_1-a_2+2(k+1)\eta,a_1-a_2+2(2k-N)\eta)}{\theta(2(k-N)\eta,a_1-a_2+2(k-N)\eta,a_1-a_2+2(2k+2-N)\eta)}
\,\Gamma_k.$$
Upon iteration, this yields
$$\Gamma_k=\frac{[\lambda]}{[\lambda+2k]}\frac{[1,\lambda+N+1]_k}{[-N,\lambda]_k}\,\Gamma_0, $$
where $\lambda=(a_1-a_2-2N\eta)/2\eta$.

We are thus reduced to computing $\Gamma_0=C^{-1}\langle e_0,f_0\rangle$, for
which we observe that 
$$f_0=e^{\pi iN(\tau-1)/2}K_{a_1-\frac12 N\eta},$$
and thus, by Theorem \ref{sxc},
$$\Gamma_0=e^{\pi iN(\tau+1)/2}\,e_0(a_1-\frac 12 N\eta)=
e^{\pi iN(\tau-1)/2}\left[\lambda+1,(a_1+a_2-N\eta)/{2\eta}\right]_N.$$
This completes the proof of Theorem \ref{bbt}.

\section{Elliptic hypergeometric series}
\label{ehs}

By the general theory of reproducing kernel Hilbert spaces (which
reduces to linear algebra in the present,
finite-dimensional, case) we have, in the notation of Theorems
\ref{sxc} and \ref{bbt},
\begin{equation}\label{rke}K_v(u)
=\sum_{k=0}^N\frac{e_k(u)\overline{f_k(v)}}{\Gamma_k}.
\end{equation}
Writing this out explicitly, using \eqref{fc}, one
obtains after simplification
\begin{multline}\label{js}
\sum_{k=0}^N\frac{[a+2k]}{[a]}\frac{[a,b,c,d,e,-N]_k}
{[1,a+1-b,a+1-c,a+1-d,a+1-e,a+1+N]_k}\\
= 
\frac{[a+1,a+1-b-c,a+1-b-d,a+1-c-d]_N}{[a+1-b,a+1-c,a+1-d,a+1-b-c-d]_N},
\end{multline}
where
\begin{multline*}(a,b,c,d,e)=\frac 1{2\eta}\left(a_1-a_2-2N\eta,-a_2-\frac 12N\eta+\eta+\bar
v+\frac12+\frac\tau2,\right.\\
\left.-a_2-\frac 12N\eta+\eta-\bar
v-\frac12-\frac\tau2, a_1-\frac12N\eta+u,a_1-\frac12N\eta-u\right).\end{multline*}
These parameters satisfy
 $b+c+d+e=2a+N+1$ but are otherwise generic. The summation formula \eqref{js} was first obtained
 by Frenkel and Turaev \cite{ft}. When $p=0$ it reduces to the
 Jackson sum \cite[(II.22)]{gr}.
Its appearance here is another example of the connection between
elliptic quantum groups and elliptic hypergeometric series, see further \cite{knr,r}.
Note also that, conversely, taking the
scalar product of both sides of \eqref{rke} with $f_l$ gives
$$f_l(v)=\sum_{k=0}^N\frac{\langle
  f_l,e_k\rangle{f_k(v)}}{C\bar\Gamma_k} $$
and  we recover $\langle e_k,f_l\rangle=\delta_{kl}C\Gamma_k$. Thus, if 
we are willing to assume \eqref{js}, we  obtain an alternative proof of
Theorem
  \ref{bbt} as a consequence of Theorem \ref{sxc}.

\section{Elliptic hypergeometric integrals}
\label{eis}

Theorem \ref{sxc} is equivalent to the integral formula
\begin{equation*}\iint_{\mathbb C/(\mathbb Z+\tau\mathbb
    Z)} K(v,\bar z)K(z,\bar w)M(z,\bar z)\,dxdy=CK(v,\bar w),\qquad
  z=x+iy. \end{equation*}
We write this out explicitly, using the notation
$$\Gamma(x;p,q)=\prod_{j,k=0}^\infty\frac{1-p^{j+1}q^{k+1}/x}{1-p^jq^kx}$$
for Ruijsenaars' elliptic gamma function \cite{ru}, which satisfies
$$\frac{\Gamma(qx;p,q)}{\Gamma(x;p,q)}=(x,p/x;p)_\infty,\qquad
\Gamma(x;p,q)\Gamma(pq/x;p,q)=1. $$
After making the substitution 
$e^{2\pi iz}\mapsto z$ and similarly for $v$ and $\bar w$, we obtain after simplification
the identity
\begin{multline}\label{bdi}
\iint_{p<|z|<1}\frac{\Gamma(-p^\frac12q^{\frac{N+1}{2}}v^\pm\bar
  z^\pm,- p^\frac12 q^{\frac{N+1}{2}}\bar w^\pm z^\pm;p,q)}
{\Gamma(z^2,\bar z^2,pz^{-2},p\bar z^{-2},- p^\frac 12q^{\frac{N+3}{2}}
  z^\pm \bar z^\pm;p,q)}\,\frac{dxdy}{|z|^4}\\
=\frac{2\pi\log(q)
  p^{-1/2}q^{(N+1)/2}}{(p,p,q^{N+1},pq^{-N-1};p)_\infty}\,\Gamma(-p^\frac12q^{\frac{N+1}{2}}v^\pm \bar w^\pm;p,q),\qquad z=x+iy,
\end{multline}
in  standard short-hand notation analogous to \eqref{tsh}.
We have proved this for $\eta\in\mathbb R_N\cup \mathbb I_N$, but it
extends immediately to the region
$$|\operatorname{Re}(\eta)|<1/2(N+1),\qquad
|\operatorname{Im}(\eta)|<\tau/2(N+1)i$$
 or, equivalently,
$$p^{-1/(N+1)}<|q|<p^{1/(N+1)},\qquad|\operatorname{arg}(q)|<2\pi/(N+1).$$

The identity \eqref{bdi}
 is quite similar in structure to  Spiridonov's elliptic beta
integral \cite{sp},
 and to its double and multiple
extensions conjectured by van Diejen and Spiridonov \cite{ds} and proved
by Rains \cite{rai1}. Initially, we tried to prove Theorem~\ref{sxc} by
deducing \eqref{bdi} from such known results, but we did not succeed
with this approach.

\section{Elliptic $6j$-symbols}
\label{esjs}

In \cite{r}, we studied  the change of base
coefficients $R_k^l=R_k^l(a,b,c,d;N;q,p)$ occurring in
\begin{equation}\label{rex}e_k^N(x;a,b)=\sum_{l=0}^NR_k^l\,
  e_l^N(x;c,d). \end{equation} 
It turned out that they can be identified with analytically continued
elliptic $6j$-symbols.  Note that, because of the
different normalizations used for theta functions, the quantity denoted 
$R_k^l(a,b,c,d;N;q,p)$ here equals
$$q^{(l-k)(k+l-N)}e^{2\pi i(-ak-b(N-k)+cl+d(N-l))}
R_k^l(e^{2\pi ia},e^{2\pi ib},e^{2\pi ic},e^{2\pi id};N;q,p)
 $$
in the notation of \cite{r}. 

Using Theorem \ref{bbt}, we may interpret 
$R_k^l$  as a scalar product of two basis
vectors. Namely, taking the scalar product of both sides of
\eqref{rex}
with 
$$\sigma
e_l^N(x;\eta(1-N)\mp\bar d,\eta(1-N)\mp\bar c)$$
 gives
\begin{multline}\label{rsp}R_k^l(a,b,c,d;N;q,p)=\frac{C^{-1}e^{\pi i N(1-\tau)/2}}
{[\lambda  +1,(c+d)/2\eta]_N}
\frac{[\lambda+2l]}{[\lambda]}\frac{[-N,\lambda]_l}{[1,\lambda+N+1]_l}\\
\times\left\langle e_k^N(x;a,b) ,\sigma
e_l^N(x;\eta(1-N)\mp\bar d,\eta(1-N)\mp\bar c)\right\rangle,
\end{multline}
where $\lambda=(c-d-2N\eta)/2\eta$. 
As before, the minus sign is taken for $\eta\in\mathbb R_N$ and the
plus sign for $\eta\in\mathbb I_N$.

The main interest in \eqref{rsp} is that it explains the 
self-duality of elliptic $6j$-symbols.  To this end, we apply the symmetry
$\langle f,g\rangle=\overline{\langle g,f\rangle}$
to \eqref{rsp}, using that
$\sigma^\ast=\sigma$ and that, by \eqref{ec},
$$\overline{R_k^l(a,b,c,d;N;q,p)}=R_k^l(\pm\bar a,\pm\bar b,\pm\bar c,\pm\bar d;N;q,p).$$
Combining these facts we obtain the end result
\begin{multline*}
R_k^l(a,b,c,d;N;q,p)\\
=\frac{[\mu  +1,(a+b)/2\eta]_N}{[\lambda  +1,(c+d)/2\eta]_N}
\frac{[\lambda+2l]}{[\lambda]}\frac{[-N,\lambda]_l}{[1,\lambda+N+1]_l}
\frac{[\mu]}{[\mu+2k]}\frac{[1,\mu+N+1]_k}{[-N,\mu]_k}\\
\times R_l^k(\eta(1-N)-d,\eta(1-N)-c,\eta(1-N)-b,\eta(1-N)-a;N;q,p), 
\end{multline*}
where $\lambda=(c-d-2N\eta)/2\eta$ and
$\mu=(a-b-2N\eta)/2\eta$.
This symmetry 
 is not obvious from \eqref{rex}, although it is
clear from  the explicit expression for $R_k^l$ as an elliptic
hypergeometric series given in \cite[Theorem 3.3]{r}.

\section*{Appendix. Sklyanin's generators}
\setcounter{section}{1}
\renewcommand{\thesection}{\Alph{section}}
\setcounter{equation}{0}

In this appendix we provide the details of the relation between
the operators \eqref{rso} and Sklyanin's generators. This is mostly
 based on personal communication from
Eric Rains.

Sklyanin  operators \cite[Theorem 2]{sk2} have the form
$$S_if(x)=\frac{s_i(x-\frac12N\eta)f(x+\eta)-s_i(-x-\frac12N\eta)f(x-\eta)}
{\theta(2x)},$$
where, in our notation,
\begin{equation*}\begin{split}s_0(x)&=\theta(\eta,2x),\\
s_1(x)&=\theta\left(\eta+\frac12,2x+\frac12\right),\\
s_2(x)&=e^{\pi i(\frac12+\frac\tau2+\eta+2x)}
\theta\left(\eta+\frac12+\frac\tau2,2x+\frac12+\frac\tau2\right),\\
s_3(x)&=-e^{\pi i(\frac\tau2+\eta+2x)}
\theta\left(\eta+\frac\tau2,2x+\frac\tau2\right).
\end{split}\end{equation*}
Using \eqref{dupl}, it is easy to
check that 
\begin{equation}\label{sg}\begin{split}
S_0&=\frac{ip^{1/8}\theta(\eta)}{(p;p)_\infty^3}\,\Delta\left(0,\frac12,\frac\tau2,-\frac12-\frac\tau2\right),\\
S_1&=-\frac{ip^{1/8}\theta(\eta+\frac
  12)}{(p;p)_\infty^3}\,\Delta\left(\frac 14,-\frac
14,\frac 14+\frac\tau 2,-\frac 14-\frac \tau 2\right),\\
S_2&=\frac{ip^{1/8}e^{\pi i\eta}\theta(\eta+\frac
  12+\frac\tau2)}{(p;p)_\infty^3}\,\Delta\left(\frac 14+\frac\tau 4,\frac 14-\frac \tau 4,-\frac 14+\frac\tau 4,-\frac 14-\frac \tau 4\right),\\
S_3&=\frac{ip^{1/8}e^{\pi
    i\eta}\theta(\eta+\frac\tau2)}{(p;p)_\infty^3}\,\Delta\left(\frac\tau4,-\frac\tau4,\frac\tau4+\frac12,-\frac\tau4-\frac12\right).
\end{split}\end{equation}

Conversely,  
every  $\Delta(a)$ is a linear combination of the $S_i$. Namely,
\begin{equation*}\begin{split}
\Delta(a)&=\frac 12\left\{\frac{\theta(a_1+a_4,a_2+a_4,a_3+a_4)}{\theta(\eta)}\,S_0\right.\\
&\qquad-\frac{\theta(a_1+a_4+\frac 12,a_2+a_4+\frac 12,a_3+a_4+\frac 12)}{\theta(\eta+\frac 12)}\,S_1\\
&\qquad-e^{\pi i(\frac\tau2+\frac12+2a_4-\eta)}\frac{\theta(a_1+a_4+\frac12+\frac\tau2,a_2+a_4+\frac12+\frac\tau2,a_3+a_4+\frac12+\frac\tau2)}{\theta(\eta+\frac12+\frac\tau2)}\,S_2\\
&\left.\qquad+e^{\pi i(\frac\tau2+2a_4-\eta)}\frac{\theta(a_1+a_4+\frac\tau2,a_2+a_4+\frac\tau2,a_3+a_4+\frac\tau2)}{\theta(\eta+\frac\tau2)}\,S_3\right\}.
\end{split}\end{equation*}
Writing this out explicitly, one is reduced to  the  theta
function identity 
\begin{multline*}
\theta(x+a_1,x+a_2,x+a_3,x+a_4)\\
\begin{split}&=\frac 12\Big\{\theta(a_1+a_4,a_2+a_4,a_3+a_4,2x)\\
&\qquad-\theta(a_1+a_4+\frac 12,a_2+a_4+\frac 12,a_3+a_4+\frac
12,2x+\frac12)\\
&+e^{\pi i(\tau+2a_4+2x)}\theta(a_1+a_4+\frac{1+\tau}2,a_2+a_4+\frac{1+\tau}2,a_3+a_4+\frac{1+\tau}2,2x+\frac{1+\tau}2)\\
&\qquad-e^{\pi i(\tau+2a_4+2x)}\theta(a_1+a_4+\frac\tau2,a_2+a_4+\frac\tau2,a_3+a_4+\frac\tau2,2x+\frac\tau2)\Big\},\end{split}
\end{multline*}
which is obtained from \eqref{ji} after substituting
$b=(a_1+a_4,a_2+a_4,a_3+a_4,-2x)$.


\begin{thebibliography}{99} 



\bibitem[DJKMO]{d2}  E.\ Date, M.\ Jimbo, A.\ Kuniba, T.\ Miwa and M.\
Okado, \emph{Exactly solvable SOS models II. Proof of the
star-triangle  relation and combinatorial identities}, in 
M.\ Jimbo et al.\ (eds.), Conformal
Field  Theory and Solvable
  Lattice Models,  17--122, Adv.\ Stud.\ Pure Math.\ 16, Academic
Press,  Boston, MA, 1988. 

\bibitem[DJMO]{d3} E.\ Date, M.\ Jimbo, T.\ Miwa and M.\ Okado, 
\emph{Fusion of the eight vertex SOS model}, Lett.\ Math.\ Phys.\ 12
(1986),  209--215; \emph{Erratum and addendum},  
Lett.\ Math.\ Phys.\ 14 (1987),  97.


\bibitem[DS]{ds} J.\ F.\  van Diejen and V.\ P.\ Spiridonov,
  \emph{Elliptic Selberg integrals},  Internat.\ Math.\ Res.\ Notices
  20 (2001),   1083--1110. 

\bibitem[FT]{ft} I.\ B.\  Frenkel and V.\ G.\ Turaev, {\it Elliptic
solutions of the Yang--Baxter equation and modular hypergeometric
functions}, in V.\ I.\ Arnold et al.\ (eds.), The
   Arnold--Gelfand Mathematical Seminars, 171--204, 
Birkh\"auser, Boston, 1997.

\bibitem[GR]{gr} G.\ Gasper and M.\ Rahman,  Basic Hypergeometric Series,
Cambridge University Press, Cambridge, 1990.


\bibitem[KNR]{knr} E.\ Koelink, Y.\ van Norden and H.\ Rosengren,
\emph{Elliptic $\mathrm{U}(2)$ quantum group and elliptic
  hypergeometric series}, Comm.\ Math.\ Phys.\ 245 (2004), 519--537. 


\bibitem[M]{m} D.\ Mumford, 
Tata Lectures on Theta, Vol.\ I., 
Birkh\"auser,  Boston,  1983.

\bibitem[R1]{rai1} E.\ M.\ Rains, \emph{Transformations of elliptic
hypergeometric integrals}, math.QA/0309252.

\bibitem[R2]{rai} E.\ M.\ Rains, \emph{$BC_n$-symmetric abelian
functions}, math.CO/0402113. 

\bibitem[Ro1]{r2} H.\ Rosengren, \emph{Elliptic hypergeometric series
on root systems}, Adv.\ Math.\ 181 (2004), 417--447. 

\bibitem[Ro2]{r} H.\ Rosengren, \emph{An elementary approach to
    $6j$-symbols (classical, quantum, rational, trigonometric, and
    elliptic)},  math.CA/0312310.

\bibitem[Ru]{ru} S.\ N.\ M.\ 
Ruijsenaars, 
\emph{First order analytic difference equations and integrable quantum
  systems}, 
J.\ Math.\ Phys.\ 38 (1997),  1069--1146.


\bibitem[S1]{sk1} E.\ K.\ Sklyanin,  
\emph{Some algebraic structures connected with the Yang--Baxter
equation},  Functional Anal.\ Appl.\ 16 (1982),  263--270. 


\bibitem[S2]{sk2} E.\ K.\ Sklyanin,  
\emph{Some algebraic structures connected with the Yang--Baxter
equation. Representations of quantum algebras},  Functional
Anal.\ Appl.\ 17 (1983),  273--284. 
 
\bibitem[Sp]{sp}  V.~P.~Spiridonov, 
\emph{On the elliptic beta function}, Russian Math.\ Surveys 56 (2001),
185--186. 

\bibitem[WW]{ww} E.\ T.\ Whittaker and G.\ N.\ Watson, 
A Course on Modern Analysis, 4th ed., Cambridge University 
Press, Cambridge, 1927 (reprinted 1996).



\vskip 3mm
\end{thebibliography}
\end{document}